\documentclass[12pt]{article}
\usepackage{amsmath,amsfonts,amssymb,color,a4,epsfig,graphics}
\usepackage[latin1]{inputenc}
\usepackage[english,frenchb]{babel}

\newcommand{\R}{{\mathbb{R}}}

\def\ESA{~essentiellement auto-adjoint~}
\def\ssi{~si et seulement si~}

\def\preuve{\begin{proof}}
\newtheorem{defi}{D\'efinition}[section]

\newtheorem{lemm}{Lemme}[section]
\newtheorem{prop}{Proposition}[section]
\newtheorem{rem}{Remarque}[section]
\newtheorem{crit}{Crit\`ere}[section]

\newtheorem{theo}{Th\'eor\`eme}[section]

\newenvironment{absf}{\noindent {\it \textbf{R\'esum\'e:}}}

\newenvironment{demo}{\noindent {\it Preuve.--}
      \begin{quotation}\noindent}{\end{quotation}\hfill$\square $}
\newenvironment{demorefpos}{\noindent {\it Preuve du lemme \ref{pos}}
      \begin{quotation}\noindent}{\end{quotation}\hfill$\square $}
\newenvironment{demorefharnak}{\noindent {\it Preuve du lemme \ref{harnak}}
      \begin{quotation}\noindent}{\end{quotation}\hfill$\square $}
\newenvironment{demorefharmo}{\noindent {\it Preuve du th\'eor\`eme \ref{harmo}}
      \begin{quotation}\noindent}{\end{quotation}\hfill$\square $}
\newenvironment{demorefunit}{\noindent {\it Preuve du th\'eor\`eme \ref{unit}}
      \begin{quotation}\noindent}{\end{quotation}\hfill$\square $}

\newenvironment{remer}{\noindent {\textbf{Remerciements:}}}

\begin{document}

\title{Laplaciens de graphes infinis I \footnote{ \textbf{Math Subject Calssification (2000):} 05C63, 05C50, 05C12, 35J10, 47B25.}\\
Graphes m\'etriquement complets\footnote {\textbf{Mots clés:} graphe infini,
Laplacien de graphe pond\'er\'e, op\'rateur de Schr\"{o}dinger, essentiellement
auto-adjoint.}}
\author{Nabila Torki-Hamza\footnote{%
Facult\'e des Sciences de Bizerte, Universit\'e 7 Novembre \`a
Carthage,(Tunisie); \texttt{nabila.torki-hamza@fsb.rnu.tn}; \texttt{%
torki@fourier.ujf-grenoble.fr}}}

\maketitle

\begin{center}
\small{\emph{ A la m\'emoire de mon p\`ere Pr. Dr. Ing. B\`echir Torki (1931-2009)}}
\end{center}

\vspace{0.5cm}

\begin{absf}
\small{\emph{ On introduit le Laplacien $\Delta_{\omega,c}$ d'un graphe $G$ localement fini pond\'er\'e \`a la fois sur
les sommets et sur les ar\^etes, ainsi que la notion d'op\'erateur de Schr\"{o}dinger $\Delta_{1,a}+W$.
Pour les graphes \`a poids constants sur les sommets, on \'etend
un r\'esultat de Wojciechowski pour le Laplacien et un r\'esultat
de Dodziuk pour les op\'erateurs de Schr\"{o}dinger concernant le caract\`ere
essentiellement auto-adjoint.}}
\newline
\small{\emph{Le r\'esultat principal de ce travail \'etablit  que pour les
graphes pond\'er\'es \`a valence born\'ee  et m\'etriquement
complets, le Laplacien  $\Delta_{\omega,c}$ est
 essentiellement auto-adjoint, et il en va de m\^eme pour l' op\'erateur
 $\Delta_{1,a}+W$ pourvu que la forme quadratique associ\'ee soit minor\'ee.
La preuve fait appel \`a la construction d'une fonction harmonique strictement
positive qui permet d'\'ecrire l'op\'erateur de Schr\"{o}dinger $\Delta_{1,a}+W$
comme un Laplacien \`a poids  $\Delta_{\omega,c}$ \`a transformation unitaire pr\`es.}}
 \end{absf}

\maketitle

\section{Introduction}

Cet article est le premier d'une s\'erie de 3 articles  (les deux autres sont
 \cite {C-To-Tr-1} et \cite {C-To-Tr-2}) qui sont consacr\'es \`a la th\'eorie spectrale des
op\'erateurs de type Laplacien et Schr\"{o}dinger sur les graphes
infinis. Nous \'etendons au cas des graphes infinis un certain
nombre de r\'esultats classiques sur les Laplaciens et op\'erateurs de
Schr\"{o}dinger sur les vari\'et\'es Riemanniennes non compactes.
\newline
Un des r\'esultats principaux de cet article, le th\'eor\`eme \ref{lapcompl},
est que le Laplacien d'un graphe pond\'er\'e \`a valence born\'ee
m\'etriquement complet est essentiellement auto-adjoint.
Le th\'eor\`eme 1.3.1 de \cite{Wo} et le th\'eor\`eme 3.1 de \cite{Jo}
en sont des cas particuliers.
La notion de compl\'etude utilis\'ee pour les graphes est relative
\`a une distance fabriqu\'ee \`a l'aide des poids sur les sommets et sur les ar\^etes.
\newline

Dans le 2\`eme article, nous nous int\'eresserons au cas
des graphes m\'etriquement non complets et donnerons des conditions
de croissance du potentiel assurant qu'un op\'erateur de Schr\"{o}dinger
est essentiellement auto-adjoint. Et le 3\`eme article
traite le cas avec champ magn\'etique.
\newline

La recherche de conditions pour qu'un op\'erateur de Schr\"{o}dinger soit \ESA
est un probl\`eme classique de la physique math\'ematique. Beaucoup de
travaux sont consacr\'es \`a l'op\'erateur de Schr\"odinger dans ${\mathbb{R}}%
^{n}$~; selon \cite{B-M-S}, le premier article sur ce sujet est
celui de Weyl \cite{Wey} et les livres \cite{RS} contiennent les
r\'esultats classiques. Plus tard Gaffney a prouv\'e dans \cite{Ga1} et \cite{Ga2} (voir
aussi \cite{Ch}, \cite{St}) que le Laplacien d'une vari\'et\'e Riemannienne
compl\`ete est \ESA. Et dans \cite{Ol}, (voir aussi \cite{Sh}, \cite{Shu}), il est prouv\'e qu'un
 op\'erateur de Schr\"odinger sur une vari\'et\'e Riemannienne compl\`ete
 est \ESA d\`es que le potentiel v\'erifie une condition de minoration.
 \newline

Plusieurs d\'efinitions de Laplaciens sur les graphes, analogues \`a celle
du Laplacien de Beltrami des vari\'et\'es Riemanniennes, ont \'et\'e
propos\'ees telles que les Laplaciens de graphes quantiques (voir
\cite{E-Ke-Ku-S-T}, \cite{Ku}, \cite{Ca}) et les Laplaciens combinatoires (voir \cite{Col},
\cite{Wo}, \cite{Go-Sch},\cite{Jo}), ou Laplaciens physiques(voir \cite{Web}).
\newline

Dans ce travail, un type diff\'erent de Laplacien, not\'e $\Delta_{\omega,c}$~,
est introduit pour un graphe localement fini pond\'er\'e
par un poids $\omega$ sur les sommets et une conductance $c$ sur
les ar\^etes. Il g\'en\'eralise aussi bien le ``Laplacien combinatoire''
dans \cite{Wo} (qui n'est autre que $\Delta_{1,1}$) que le ``graph Laplacian''
dans \cite{Jo} (qui est $\Delta_{1,c}$)~.
Cette notion a \'et\'e d\'ej\`a introduite dans le cas des
graphes finis \cite{To}, \cite{Col}.
\newline

On donne dans la section \ref{prel}
certaines  propri\'et\'es imm\'ediates du Laplacien $\Delta_{\omega,c}$
et on montre qu'il est  unitairement \'equivalent par une transformation
diagonale \`a un op\'erateur de Schr\"{o}dinger de la forme $\Delta_{1,a}+W$.
\newline

Dans la section \ref{prr}, en s'inspirant de la m\'ethode de \cite{Wo}, on
d\'emontre que, si le poids $\omega$ est constant, l'op\'erateur $%
\Delta_{\omega,c}$ est \ESA. La m\'ethode utilis\'ee permet aussi de prouver
que si on lui ajoute un potentiel $W$ minor\'e il reste \ESA.
\newline

La section \ref{fonctharm} est une partie consacr\'ee \`a la construction d'une
fonction $\Phi$ strictement positive et harmonique pour un op\'erateur de
Schr\"{o}dinger positif. Cette construction fait appel \`a l'in\'egalit\'e
de Harnak locale, \`a la r\'esolution d'un probl\`eme de Dirichlet
et au principe du minimum pour les graphes. Une telle fonction  $\Phi$ est utilis\'ee
dans la section \ref{operunit} pour montrer le r\'esultat important
que  \textit{tout op\'erateur de Schr\"{o}dinger positif est unitairement
\'equivalent \`a un Laplacien}.
\newline

On consid\`ere, dans la section \ref{complesa}, le cas des graphes \`a valence born\'ee.
Pour un op\'erateur de Schr\"{o}dinger donn\'e $%
\Delta_{1,a}+W$~, on introduit une distance $\delta_{a}$ sur le graphe, et on montre que
\textit{si le graphe est complet pour cette distance et si la forme
quadratique associ\'ee \`a  cet op\'erateur est born\'ee inf\'erieurement,
l'op\'erateur de Schr\"{o}dinger $\Delta_{1,a}+W$ est \ESA.}
Ce r\'esultat n'est pas un cas particulier  de celui du th\'eor\`eme \ref%
{op.eaa}, on \'etudie un contre-exemple pour cela.
On d\'eduit aussi, dans cette section, un r\'esultat analogue pour
le Laplacien $\Delta_{\omega,c}$~, il s'agit du r\'e%
sultat principal de cet article qui est une g\'en\'eralisation du th\'eor\`eme de
Gaffney au cas des graphes m\'etriquement complets.

\section{Pr\'eliminaires}\label{prel}

Soit $G$ un graphe connexe, infini et localement fini. Nous d\'esignons par $V$
l'ensemble de ses sommets et par $E$ celui de ses ar\^etes. Pour deux sommets $x$
et $y$ de $V$~, nous notons  $x\sim y$ s'ils sont reli\'es par
une ar\^ete qui sera d\'esign\'ee par $\lbrace x,y\rbrace \in E$~.
Lorsque dans certains calculs, le graphe $G$ est suppos\'e orient\'e,
nous notons $[x,y]$ l'ar\^ete d'origine $x$ et d'extr\'emit\'e $y$~,
et nous d\'esignons par $\overline{E}$ l'ensemble des ar\^etes orient\'ees.
Il est \`a signaler qu'aucun r\'esultat ne d\'epend de l'orientation.
\newline

L'ensemble des fonctions sur $V$ est not\'e par
 $$C\left( V\right) =\lbrace f:~V \longrightarrow {\mathbb{R}}\rbrace$$ et celui des
fonctions \`a support fini par $C_0\left( V\right)$. Soit $%
\omega:~V\longrightarrow {\mathbb{R}}_+^{\star}$  une fonction poids sur les
sommets, consid\'erons l'ensemble
\begin{equation*}
l_{\omega}^2 \left( V\right)= \left\lbrace f:~V\longrightarrow {\mathbb{R}}~;~
\sum_{x\in V}\omega_x^2\vert f \left( x \right) \vert^2\ <
\infty\right\rbrace ~.
\end{equation*}
L'espace $l_{\omega}^2 \left( V\right)$ muni du produit scalaire donn\'e
par:
\begin{equation*}
\langle f,g \rangle_{ l_{\omega}^2 } =\sum_{x\in V} \omega_x^2 f\left(
x\right).g\left( x\right)
\end{equation*}
est un espace de Hilbert isomorphe \`a
\begin{equation*}
l^2 \left( V\right)= \left\lbrace f:~V\longrightarrow {\mathbb{R}}~;~\sum_{x\in
V}\vert f \left( x\right) \vert^2\ < \infty\right\rbrace
\end{equation*}
par la transformation unitaire
\begin{equation*}
U_\omega :~l_{\omega}^2 \left( V\right)\longrightarrow l^2 \left( V\right)
\end{equation*}
d\'efinie par
\begin{equation*}
U_\omega\left( f\right) =\omega f~.
\end{equation*}

\begin{rem}
Si le poids $\omega$ est constant \'egal \`a $\omega_0 \ >0$ (c'est \`a dire
que: pour tout x$\in V$, on a $\omega_x=\omega_0$)~, alors
\begin{equation*}
l_{\omega_0}^2 \left( V\right)=l^2\left( V\right)~.
\end{equation*}

\end{rem}

\begin{defi}\label{laplacien}
\label{lap} Le \emph{Laplacien du graphe $G$ pond\'er\'e} par un poids
$\omega:V\longrightarrow {\mathbb{R}}_+^{\star}$ sur les sommets et
une conductance $c:E\longrightarrow {\mathbb{R}}_+^{\star}$ sur les
ar\^etes, est l' op\'erateur sur $l_{\omega}^2\left( V\right)$~, qu'on note
$\Delta_{\omega,c}$~, donn\'e par :
\begin{equation*}
\left( \Delta_{\omega,c}f\right)\left( x\right) = \dfrac{1}{\omega_{x}^{2}}
\sum_{ x\sim y } c_{\lbrace x,y\rbrace}\left( f\left(
x\right) -f\left(y \right) \right)
\end{equation*}
pour tout $f\in l_{\omega}^2\left( V\right)$ et pour tout sommet $x$ de $V$~.
\end{defi}

\begin{rem}
Voici quelques propri\'et\'es simples de ces Laplaciens dont certaines sont
inspir\'ees de \cite{Col} et \cite{Dod}~:

\begin{enumerate}
\item L'op\'erateur $\Delta_{\omega,c}$ est sym\'etrique sur  $%
l_{\omega}^2\left( V\right)$ avec domaine $C_0(V)$~; la forme quadratique
associ\'ee est
\begin{equation*}
Q_c (f)=\sum _{ {\lbrace x,y\rbrace}\in E } c_{{\lbrace x,y\rbrace}%
}(f(x)-f(y))^2~
\end{equation*}
est positive.

\item L'op\'erateur $\Delta_{\omega,c}$ s'annule pour les fonctions
constantes, d\`es que le poids $\omega$ appartient \`a $l^2\left( V\right)$~.

\item Le graphe $G$ \'etant localement fini,  cet op\'erateur est bien
d\'efini sur $C_0(V)$~, car les sommes qui interviennent sont finies.

\item C'est un op\'erateur local, au sens que $\left(\Delta_{\omega,c}f
\right)\left( x\right)$ ne d\'epend que des valeurs de $f$ aux sommets
voisins de $x$.  On peut ainsi consid\'erer le Laplacien $\Delta_{\omega,c}$
comme un \textrm{op\'erateur diff\'erentiel} sur le graphe $G$~.

\item Cet op\'erateur est \textrm{elliptique,} puisque pour chaque ar\^ete $%
\lbrace x,y\rbrace$ de $G$~, le coefficient $c_{\lbrace x,y\rbrace}$ n'est
pas nul.

\item La fonction $c$ ne d\'epend pas de l'orientation de l'ar\^ete, et on
a: $c_{\lbrace x,y\rbrace}=c_{\lbrace y ,x\rbrace}$~, pour tous sommets $x$
et $y$~.
\end{enumerate}
\end{rem}

Pour se ramener au m\^eme espace de fonctions $l^2\left(
V\right)$~, nous utilisons la transformation unitaire
$U_\omega$~. Et plus pr\'ecis\'ement, la proposition
\ref{lap-unit}, affirme que $\Delta_{\omega,c}$ est unitairement
\'equivalent \`a un op\'erateur de Schr\"{o}dinger du graphe $G$
dont ci-dessous la d\'efinition .

\begin{defi}
Un \emph{op\'erateur de Schr\"odinger du graphe} $G$ est un op\'erateur
de la forme $\Delta_{1,a} +W $ op\'erant sur $l^2\left( V\right)$~,
o\`u $W$ est une fonction r\'eelle sur $V$
et $a$ est une fonction strictement positive sur $E$~.
\end{defi}

\begin{prop}
\label{lap-unit}  Si
\begin{equation*}
\widehat{\Delta}= U_\omega \Delta_{\omega,c}U_\omega ^{-1}~,
\end{equation*}
alors $\widehat{\Delta}$ est un op\'erateur de Schr\"odinger
de $G$ et on a plus pr\'ecis\'ement
\begin{equation*}
\widehat{\Delta}=\Delta_{1,a}+W
\end{equation*}
o\`u $a$ est la fonction strictement positive sur $E$ donn\'ee par:
\begin{equation*}
a_{\lbrace x,y\rbrace}=\dfrac{c_{{\lbrace x,y\rbrace}}}{\omega_x\omega_y}
\end{equation*}
et le potentiel $W:V\longrightarrow {\mathbb{R}}$ est donn\'e par~:
\begin{equation*}
W=-\dfrac{1}{\omega} \Delta_{1,a}\omega~.
\end{equation*}
\end{prop}

\begin{demo}
Pour $g\in\ C_0\left( V\right)$, calculons $\left( \widehat{\Delta}g\right)
\left(x \right) $~.
\begin{align}
\left( \widehat{\Delta}g\right) \left(x \right)&= \omega_{x}\left(
\Delta_{\omega,c}\mathit{U_\omega}^{-1}g\right) \left( x\right)  \notag \\
&=\dfrac{1}{\omega_{x}}\sum_{\lbrace x,y\rbrace\in E}c_{\lbrace
x,y\rbrace}\left( \dfrac{g\left(x \right) } {\omega_{x}}-\dfrac{g\left(y
\right) }{\omega_{y}}\right)  \notag \\
&=\sum_{\lbrace x,y\rbrace\in E} \dfrac{c_{{\lbrace x,y\rbrace}}} {%
\omega_{x}\omega_{y}}\left( g\left( x\right)-g\left( y\right)\right)
+g\left( x\right)\dfrac{1}{\omega_{x}}\sum_{\lbrace x,y\rbrace\in E}
c_{\lbrace x,y\rbrace}\left(\dfrac{1} {\omega_{x}}-\dfrac{1}{\omega_{y}}
\right)  \notag \\
&=\left( \Delta_{1,a} g\right)\left( x\right)+W\left(x \right)g\left(
x\right)~.  \notag
\end{align}
o\`u $\Delta_{1,a}$ d\'esigne le Laplacien sur $G$ pond\'er\'e par la
fonction constante $\omega\equiv1$ sur $V$ et par
la fonction strictement positive $a$ sur $E$ donn\'ee par:
\begin{equation*}
a_{\lbrace x,y\rbrace}=\dfrac{c_{\lbrace x,y\rbrace}}{\omega_{x}\omega_{y}}
\end{equation*}
et o\`u le potentiel $W:V\longrightarrow {\mathbb{R}}$ est d\'efini par :
\begin{equation*}
W(x)=\dfrac{1}{\omega_{x}}\sum_{\lbrace x,y\rbrace\in E} c_{\lbrace
x,y\rbrace}\left(\dfrac{1}{\omega_{x}} -\dfrac{1}{\omega_{y}}\right)= -%
\dfrac{1}{\omega_x}\left( \Delta_{1,a}\omega\right) \left( x\right)~.
\end{equation*}
\end{demo}

\begin{rem}
Dans le lemme \ref{lap-unit}, il est possible d'obtenir la fonction $W$
strictement n\'egative alors que le Laplacien est positif~: on prend par exemple le
graphe $G$ avec $V={\mathbb{N}}^{\star}$ et $n\sim n+1$ pour  tout $n\in {%
\mathbb{N}}^{\star}$, et on suppose que $G$ est pond\'er\'e par le poids $\omega_{n}=\dfrac{1%
}{n}$ sur les sommets et par la conductance $c_{\lbrace n,n+1\rbrace}=\left( n+1\right)^{2}$
sur les ar\^etes; alors $W\left( n\right)=-n\left( 2n+1\right)\ <0 ~.$
\end{rem}

\section{Extension des r\'esultats de Wojciechowski et Dodziuk} \label{prr}

Nos deux premiers th\'eor\`emes sont des extensions des r\'esultats
de R.K. Wojciechowski et J. Dodziuk concernant la propri\'et\'e d'\^etre
essentiellement auto-adjoint. Rappelons d'abord cette d\'efinition.

\begin{defi}
Un op\'erateur lin\'eaire sym\'etrique non born\'e dans un espace de Hilbert
est dit essentiellement auto-adjoint s'il poss\`ede une unique
extension auto-adjointe.
\end{defi}
On \'ecrira ESA (essentially self-adjoint)~, comme abr\'eviation.
\newline
Pour d\'emontrer cette propri\'et\'e ESA, le crit\`ere \ref{ESA}~,
extrait du th\'eor\`eme X.26 dans \cite{RS}~, est tr\`es pratique.

\begin{crit}
\label{ESA} L'op\'erateur sym\'etrique d\'efini positif $%
\Delta:~C_0\left(V \right)\longrightarrow l^2( V ) $ est essentiellement
auto-adjoint si et seulement si $Ker\left(
\Delta^{\star}+1\right)=\lbrace0\rbrace$~.
\end{crit}
De la d\'efinition de l'adjoint $\Delta^{\star}$ d'un op\'erateur  $\Delta:C_0\left(V
\right)\longrightarrow l^2\left( V\right)$~, nous pouvons d\'eduire que:
\begin{equation*}
Dom\left( \Delta^{\star}\right)=  \lbrace f\in l^2\left( V\right)~;~\Delta
f\in l^2\left( V\right)\rbrace~.
\end{equation*}
Nous allons alors montrer, en utilisant une id\'ee dans la preuve du th\'eor\`eme 1.3.1
\cite{Wo}~, le r\'esultat suivant.

\begin{theo}
\label{lapeaa}  Si le poids $\omega$ est constant sur $V$ alors pour toute
conductance c sur $E$, le Laplacien $\Delta_{\omega,c}$~,
avec comme domaine $C_0\left(V \right)$~, est essentiellement
auto-adjoint.
\end{theo}

\begin{demo}
Soit $\omega_0$ un r\'eel strictement positif, et $\omega\equiv\omega_0$ sur
$V$. On consid\`ere $g$ une fonction sur $V$ v\'erifiant: $%
\Delta_{\omega_0,c} g+g=0$.\newline
Supposons qu'il existe $x_0$ dans $V$ tel que $g\left(x_0 \right)\ >0$.
\newline
L'\'egalit\'e
\begin{equation*}
\Delta_{\omega_0,c} g\left( x_0\right) +g\left( x_0\right) =0
\end{equation*}
entraine
\begin{equation*}
\dfrac{1}{\omega_0^2}\sum_{ y\sim x_0}c_{\lbrace x_0,y\rbrace}\left( g\left(
x_0\right) -g \left( y\right) \right)+g\left( x_0\right) =0~.
\end{equation*}
Donc il existe au moins un sommet $x_{1}$ pour lequel $g\left( x_0\right) \
< g\left( x_{1}\right)~,$ puisque $\omega_0\ >0$ et $c_{\lbrace x,y\rbrace} \
> 0$ pour tout $\lbrace x,y\rbrace\in E~.$ On r\'eit\`ere ensuite avec $x_{1}~,$
et on construit ainsi une suite strictement croissante de r\'eels
strictement positifs $\left( g\left(x_{n} \right) \right)_{n}~.$ Ce qui
entraine que la fonction $g$ n'est pas dans $l^2\left( V\right)~.$ Un m\^eme
raisonnement est utilis\'e en supposant $g\left(x_0\right) \ < 0~.$
\end{demo}
\begin{rem}
Le th\'eor\`eme 1.3.1 dans \cite{Wo} est un cas particulier du th\'eor\`eme \ref{lapeaa}~.
\end{rem}
On peut, avec un raisonnement analogue, d\'emontrer le th\'eor\`eme suivant.
\begin{theo}
\label{op.eaa}  Si $W:V\longrightarrow {\mathbb{R}}$ est un potentiel
minor\'e et si $\omega_0$ est un poids constant sur $V$, alors pour toute
conductance $c:E\longrightarrow {\mathbb{R_{+}}}$, l'op\'erateur de Schr\"odinger $\Delta_{\omega_0,c}+W$,
avec pour domaine $C_0\left(V\right)$, est \ESA.
\end{theo}

\begin{demo}
Soit $\kappa$ un r\'eel minorant le potentiel $W$~. On proc\`ede  comme
dans la preuve du th\'eor\`eme \ref{lapeaa}~, en consid\'erant une fonction $g$ sur $V$ v\'erifiant: $%
\Delta_{\omega_0,c} g+Wg+\kappa_{1}g=0~,$ avec $\kappa +\kappa_{1} \geq 1~.$
\end{demo}

\begin{rem}

J. Dodziuk affirme par le th\'eor\`eme 1.2 dans \cite{Dod} que $A+W$ est
essentiellement auto-adjoint si l'op\'erateur $A$ est sym\'etrique
positif born\'e sur $l^2(V)$ et si $W$ est minor\'e.\newline
Dans le th\'eor\`eme \ref{op.eaa}~, l'op\'erateur $A$ est not\'e $\Delta_{1,c}$ et nous pouvons
conclure que  ce th\'eor\`eme est plus g\'en\'eral que celui de Dodziuk,
puisque $\Delta_{1,a} +W$ est essentiellement auto-adjoint
si $W$ est minor\'e, m\^eme si l'op\'erateur $A=\Delta_{1,a}$ n'est pas
born\'e sur $l^2(V)$~, en prenant par exemple le graphe localement fini
\`a valence non born\'ee et $a\equiv1~.$
\end{rem}

\section{Construction d'une fonction harmonique sur les sommets}\label{fonctharm}

Nous allons construire une fonction $\Phi$ strictement positive et
harmonique sur les sommets qui sera utile dans la section \ref{operunit}~.

\begin{theo}
\label{harmo} Soit $P$ un op\'erateur de Schr\"odinger sur $l^2\left(V\right)$
tel que pour tout $f\in C_0\left( V\right)\setminus\lbrace
0\rbrace ~,$ $\langle Pf,f\rangle_{l^{2}} > 0$~. Alors il existe une fonction $\Phi$
strictement positive et $P-$harmonique sur $V$~.
\end{theo}

La preuve du th\'eor\`eme \ref{harmo} s'appuie sur le lemme \ref{harnak} qui
donne l'in\'egalit\'e de Harnack locale pour les graphes. Nous pr\'esentons
d'abord les d\'efinitions suivantes:

\begin{defi}
Un sous-graphe $G^{\prime}$ de $G$ est un graphe dont l'ensemble des sommets
est inclus dans $V$ et celui des ar\^etes est un sous-ensemble de $E$~.
\end{defi}

\begin{defi}
Pour un sous-graphe $K$ de $G$~, on d\'efinit :

\begin{itemize}
\item l'int\'erieur de $K$ qu'on note $\overset{\circ}{K}$
\begin{equation*}
\overset{\circ}{K}=\{ x\in K;y\sim x\Rightarrow y\in K\}
\end{equation*}

\item le bord de $K$ qu'on note $\partial K$
\begin{equation*}
\partial K= K\setminus \overset{\circ}{K}= \left\lbrace x\in K~;\exists y\in
V\setminus K,y\sim x \right\rbrace
\end{equation*}

\item $K $ est connexe si et seulement si pour tous $x $ et $y$ de $K$ il
existe des sommets $x_1,x_2,...,x_n$~, tels que
\begin{equation*}
x_i \in K,\; x_1=x,\; x_n=y,\; \{x_i,x_{i+1}\}\in E
\end{equation*}
pour tout $1\leq i\leq n-1$~.
\end{itemize}
\end{defi}

\begin{lemm}
\label{harnak} Soit $P$ un op\'erateur de Schr\"odinger sur $l^2\left(
V\right)$~. Pour tout sous-graphe fini $K$ d'int\'erieur connexe fini, il existe
une constante $k\ >0$ telle que, pour toute fonction $\varphi:V%
\longrightarrow {\mathbb{R}} $ strictement positive sur les sommets de $K$
et v\'erifiant $\ (P\varphi \ )\upharpoonright \overset{\circ}{K} \equiv 0~,$
on a:
\begin{equation*}
\dfrac{1}{k} \leq \dfrac{\varphi\left( x\right) }{\varphi\left( y\right)}%
\leq k
\end{equation*}
pour tous $x,y\in \overset{\circ}{K}~.$
\end{lemm}

La r\'esolution du probl\`eme de Dirichlet fournie par le lemme \ref{pos}
suivant, est aussi utile pour la preuve du th\'eor\`eme \ref{harmo}~.

\begin{lemm}
\label{pos} Soient $P$ un op\'erateur de Schr\"odinger du graphe $G$
tel que pour tout $f\in C_0\left( V\right)\setminus\lbrace0\rbrace$
on ait $ \langle Pf,f \rangle_{ l^2 }>0~$.
Alors pour tout sous-graphe $K$ de $G$ d'int\'erieur
connexe fini et pour toute fonction $u:\partial K\longrightarrow {\mathbb{R}}~$, il existe une
fonction unique $f$ sur $K$ v\'erifiant les deux conditions suivantes:%
\newline
\textit{(i)} $\ (Pf\ )\upharpoonright \overset{\circ}{K}\equiv0~.$ \newline
\textit{(ii)} $f\upharpoonright \partial K\equiv u~.$\newline
De plus, si $u$ est positive et non identiquement nulle, alors $f$ est
strictement positive dans $\overset{\circ}{K}~.$
\end{lemm}

Afin de prouver la stricte positivit\'e dans le lemme \ref{pos}~, nous allons
utiliser le ``principe du minimum'' pour les graphes, donn\'e par le
lemme \ref{princ-min} \cite{Dod}~.

\begin{lemm}
\label{princ-min}  Soit $P=\Delta_{1,a}+W$ un op\'erateur de Schr\"odinger
sur le graphe $G$~, avec $W\ >0$~, et soit $K$ un sous-graphe fini de $G$
d'int\'erieur connexe. On suppose que $f$ est tel que
$\langle Pf,f\rangle_{l^{2}} \geq0$ \`a l'int\'erieur de $K$ et qu'il existe un sommet
int\'erieur $x_{0}$ tel que $f\left( x_{0}\right)$ soit minimum et
n\'egatif. Alors $f $ est constante sur $K$~.
\end{lemm}

\begin{demorefpos}
\begin{itemize}
\item Pour l'unicit\'e, on suppose l'existence de deux fonctions
$f$ et $g$ \`a supports finis dans $K$ v\'erifiant les deux
conditions du th\'eor\`eme, alors il s'ensuit que $P\left(
f-g\right)\equiv0 $ \`a l'int\'erieur de $K$~, et que
$\left(f-g\right)\upharpoonright \partial K\equiv0 $~. Ce qui
entraine que $\left\langle P\left( f-g\right),f-g\right\rangle_{l^{2}}=0$~.
On en d\'eduit, par l'hypoth\`ese faite sur $P$ que $\left(
f-g\right)$ est identiquement nulle, d'o\`u l'unicit\'e.\newline

\item L'unicit\'e implique l'existence car l'espace des fonctions sur $K$
est de dimension finie.\newline

\item Pour la stricte positivit\'e, on prend $u$ positive et non
identiquement nulle et on raisonne par l'absurde pour montrer que $f$ est
strictement positive \`a l'int\'erieur de $K$~. On suppose qu'il existe un
sommet dans $\overset{\circ}{K}$ dont l'image par $f$ est n\'egative. On
consid\`ere alors un sommet $x_{0} $ r\'ealisant le minimum de $f$ sur $%
\overset{\circ}{K}$ qui est fini et connexe. On a ainsi $f\left(
x_{0}\right)\leq0 $ et $Pf$ nulle sur $\overset{\circ}{K}$~. Et d'apr\`es le
lemme \ref{harmo}~, l'application $f$ est constante n\'egative sur $K$~. Ce
qui est impossible, vu que $f\upharpoonright\partial K\equiv u$ et qu'on a
suppos\'e $u\geq0 $ et non identiquement nulle. D'o\`u $f$ est strictement
positive sur $K~.$
\end{itemize}
\end{demorefpos}

La preuve de l'in\'egalit\'e de Harnak est inspir\'ee des
d\'emonstrations du lemme 1.6 et du corollaire 2.3 dans
\cite{Dod}~, en remarquant que la constante $k$ qu'on trouve ne
d\'epend pas de la fonction $\varphi~.$

\begin{demorefharnak}
On consid\`ere un sous-graphe $K$ fini d'int\'erieur connexe, et une
fonction $\varphi :V\longrightarrow {\mathbb{R}} $ strictement positive sur
les sommets de $K$ et $P-$harmonique sur $\overset{\circ}{K}$~. Soient $x$ et
$y$ deux sommets de $\overset{\circ}{K}~.$

\begin{enumerate}
\item[(i)] On suppose d'abord que $\lbrace x,y\rbrace$ est une ar\^ete.%
\newline
Comme $\left(P\varphi \right)\left( x\right)=0$~, c'est \`a dire que
\begin{equation*}
\sum_{z\sim x}a_{\lbrace x,z\rbrace}\left[\varphi\left(
x\right)-\varphi\left( z \right) \right]  +W\left(x \right) \varphi\left(
x\right) =0~,
\end{equation*}
alors
\begin{equation*}
\left( \sum_{z\sim x}a_{\lbrace x,z\rbrace}\right)\varphi\left( x\right)
+W\left(x \right) \varphi\left( x\right) =\sum_{z\sim x}a_{\lbrace
x,z\rbrace} \varphi\left( z \right)~.
\end{equation*}
On obtient, par la positivit\'e de $\varphi$ et des $a_{\lbrace x,z\rbrace}$~,
 l'in\'egalit\'e suivante:
\begin{equation*}
\left[W\left(x \right) +\sum_{z\sim x}a_{\lbrace x,z\rbrace}\right]
\varphi\left( x\right)  \geqslant a_{\lbrace x,y\rbrace} \varphi\left(
y\right)~.
\end{equation*}
On note: ~$\alpha=\min\left\lbrace a_{\left\lbrace r,s\right\rbrace } ;r,s\in
K , r\sim s\right\rbrace$ et
\begin{equation*}
A=\sum_{r,s\in K,r\sim s} a_{\left\lbrace r,s\right\rbrace }~.
\end{equation*}
Comme $K$ est fini, on a: $\alpha \ >0$ et $A\ <\infty$~. D'o\`u en notant:
\begin{equation*}
k_{0}=\dfrac{\max \left( 0,\max_{K} W\right)+A }{\alpha}~,
\end{equation*}
on a: $k_{0}\ >0$~, et on obtient:
\begin{equation*}
\dfrac{1}{k_{0}}\leq \dfrac{\varphi\left( x\right) }{\varphi\left( y\right) }%
\leq k_{0}~.
\end{equation*}

\item[(ii)] Maintenant si les sommets $x$ et $y$ ne sont pas reli\'es par
une ar\^ete, par la connexit\'e de $\overset{\circ}{K}$~, il existe un chemin
de sommets cons\'ecutifs: $x_{1}=x,x_{2},x_{3}$,...,$x_{d}=y$ reliant $x$
\`a $y$ dans $\overset{\circ}{K}$~. On obtient alors:
\begin{equation*}
\dfrac{1}{k_{0}}\leq\dfrac{\varphi\left( x_{i}\right) } {\varphi\left(
x_{i+1}\right)}\leq k_{0}~, \quad pour 1\leq i\leq d-1~,
\end{equation*}
et par suite, en prenant $k=k_{0}^{d}$~, on obtient
\begin{equation*}
\dfrac{1}{k}\leq \dfrac{\varphi\left( x\right) }{\varphi\left( y\right) }%
\leq k~.
\end{equation*}
\end{enumerate}
\end{demorefharnak}

\begin{demorefharmo}
On suppose que $\langle Pf,f\rangle \ > 0$ pour toute fonction $f\in C_{0
}\left( V\right)\setminus\lbrace0\rbrace$~. Soit $x_{0}$ un sommet fix\'e de $%
V$, pris comme ``origine''. On consid\`ere pour $n \geq 1$, le sous-graphe $%
G_{n}$ issu de $G$, dont l'ensemble des sommets est la boule de centre $x_{0}
$ et de rayon $n$~, qu'on notera ${\mathcal{B}}_{n}~,$
\begin{equation*}
{\mathcal{B}}_{n}=\lbrace x\in V;d\left( x_{0},x\right)\leq n\rbrace
\end{equation*}
o\`u $d\left( x,y\right)$ est la distance combinatoire entre deux sommets $x$
et $y$ de $V$~, qui est le nombre d'ar\^etes du plus court chemin d'ar\^etes
reliant $x$ \`a $y$~. La boule ${\mathcal{B}}_{n}$ est connexe et on applique
le lemme \ref{pos}~, en la prenant pour $K$~, et en choisissant pour fonction $%
u$ la fonction constante \'egale \`a $1$ sur $\partial {\mathcal{B}}_{n}$~.%
\newline
On va proc\'eder en trois \'etapes:

\begin{itemize}
\item {1\`ere \'etape:} Il existe une fonction $\psi_{n}\in C_{0}\left(
V\right)$ v\'erifiant $P\psi_{n} \equiv0$~, et telle que $\psi_{n}\ >0 $ \`a
l'int\'erieur de ${\mathcal{B}}_{n }$ et constante \'egale \`a $1$ sur $%
\partial {\mathcal{B}}_{n}$~. On consid\`ere alors la fonction $\Phi_{n }\in
C_{0}\left( V\right)$ d\'efinie par:
\begin{equation*}
\Phi_{n }\left( x\right)=\dfrac{\psi_{n}\left( x\right) }{\psi_{n}\left(
x_{0}\right)~.}
\end{equation*}
Elle v\'erifie les quatre conditions suivantes :

\begin{enumerate}
\item[\textit{i}.] $\Phi_{n }\left( x_{0}\right) =1$~.

\item[\textit{ii}.] $P\Phi_{n }\equiv0$ \`a l'int\'erieur de ${\mathcal{B}}%
_{n}$~.

\item[\textit{iii}.] $\Phi_{n }\upharpoonright \partial {\mathcal{B}}%
_{n}\equiv\dfrac{1} {\psi_{n}\left( x_{0}\right) }$ \; constante strictement
positive.

\item[\textit{iv}.] $\Phi_{n} \ >0$ sur ${\mathcal{B}}_{n}~.$
\end{enumerate}

\item {2\`eme \'etape:} Soit $x$ un sommet de $V$~, on fixe $n_{0}$ tel que $x
$ soit \`a l'int\'erieur de ${\mathcal{B}}_{n_{0}}$~. Alors pour tout $n\geq
n_{0}$~, on a: ${\mathcal{B}}_{n_{0}}\subseteq {\mathcal{B}}_{n}$~. De plus $%
\Phi_{n}$ est strictement positive sur ${\mathcal{B}}_{n_{0}}$ et est $P-$%
harmonique \`a l'int\'erieur de ${\mathcal{B}}_{n_{0}}$~. Donc d'apr\`es le
lemme \ref{harnak}~, il existe une constante $k_{n_{0}}\ > 0 $ tel que l'on
ait :

\begin{equation*}
\dfrac{1}{k_{n_{0}}}\leq \dfrac{\Phi_{n}\left( x\right) } {\Phi_{n}\left(
x_{0}\right) }\leq k_{n_{0}}~.
\end{equation*}

Et comme $\Phi_{n }\left( x_{0}\right) =1$~, on obtient:

\begin{equation*}
\dfrac{1}{k_{n_{0}}}\leq \Phi_{n}\left( x\right) \leq k_{n_{0}}~.
\end{equation*}

Il s'ensuit que l'ensemble $\lbrace \Phi_{n}\left( x\right) \rbrace_{n\geq
n_{0}}$ est inclus dans le segment $\left[ \dfrac{1}{k_{n_{0}}},k_{n_{0} }%
\right]$~.

\item {\ 3\`eme \'etape:} On consid\`ere dans ${\mathbb{R}}^{V}$ le
sous-ensemble%
\begin{equation*}
C=\prod_{x\in V} \left[ \dfrac{1}{k_{n_{0}}},k_{n_{0}} \right]~.
\end{equation*}
Comme la suite $\left( \Phi_{n}\right)_{n\geq n_{0}}$ est une suite du
compact $C$~, elle admet une sous-suite convergente $\left( \Phi_{h\left(
n\right) }\right)_{n\geq n_{0}}$ pour la topologie de ${\mathbb{R}}^{V}$
vers une fonction $\Phi$ qui v\'erifie en particulier les deux conditions
suivantes:

\begin{enumerate}

\item[\textit{i}.] $\Phi$ est strictement positive sur $V$~, puisque $%
\Phi\left( x\right)\in \left[ \dfrac{1}{k_{n_{0}}},k_{n_{0}} \right]$~, pour
tout sommet $x$ de $V$~.

\item[\textit{ii}.] $P\Phi\equiv 0$ sur $V$~, puisque $\underset{%
n\rightarrow\infty}{lim} P\Phi_{h\left( n\right) }\left(x \right)=P\Phi
\left( x\right)$~, pour tout $x$ de $V$~.

\end{enumerate}
\end{itemize}
\end{demorefharmo}

La fonction $\Phi$ fournie par le th\'eor\`eme \ref{harmo}
sert \`a la construction de la transformation unitaire dans le th\'eor\`eme %
\ref{unit}~.

\section{Tout op\'erateur de Schr\"odinger positif est
unitairement \'equivalent \`a un Laplacien}\label{operunit}

On va montrer qu'un op\'erateur de Schr\"odinger, sous une condition de
positivit\'e, est unitairement \'equivalent \`a un Laplacien $%
\Delta_{\omega,c}$~.

\begin{theo}
\label{unit} Soit $P$ un op\'erateur de Schr\"odinger d'un graphe $G$~.
On suppose que $\langle Pf,f\rangle_{l^{2}} > 0$ pour toute fonction $%
f\in C_0\left( V\right)\setminus\lbrace0\rbrace $~. Alors il existe une
fonction poids $\omega:V\longrightarrow {\mathbb{R}}_{+}^{\star}$ sur $V$ et
une fonction conductance $c:E\longrightarrow {\mathbb{R}}_{+}^{\star}$
sur $E$ telles que l'op\'erateur $P$ soit unitairement \'equivalent au
Laplacien $\Delta_{\omega,c}$ sur le graphe $G$~.
\end{theo}

Pour la preuve du th\'eor\`eme \ref{unit}, on utilise une fonction $\Phi$
qui est \`a la fois strictement positive et $P-$harmonique, fournie par le
th\'eor\`eme \ref{harmo}.

\begin{demorefunit}
Consid\'erons $P =\Delta_{1,a}+W$ un op\'erateur de Schr\"odinger sur $G$
v\'erifiant les hypoth\`eses du th\'eor\`eme. Par le th\'eor\`eme \ref{harmo}~,
 il existe une fonction $\Phi$ strictement positive et $P-$harmonique sur $V
$. Alors on obtient: \;$W=-\dfrac{\Delta\Phi}{\Phi}$~.\newline
On pose $\omega=\Phi$ et pour tout $g\in l^{2}\left(V \right) $~, on pose $f=%
\dfrac{g}{\Phi}$~.\newline
On va montrer que $\left\langle Pg,g\right\rangle_{l^2}=\left\langle
\Delta_{\omega,c}f,f\right\rangle_{l^2_{\omega}} $

\begin{align}
\left\langle Pg,g\right\rangle_{l^2}&= \left\langle \Delta \left(
f\Phi\right) +Wf\Phi,f\Phi\right\rangle_{l^2}  \notag \\
&=\left\langle \Delta \left(
f\Phi\right)-f\Delta\Phi,f\Phi\right\rangle_{l^2}  \notag \\
&=\sum_{x\in V}\left[ \sum_{y\sim x} a_{\lbrace x,y\rbrace} \left(f\left(
x\right)\Phi\left( x\right)- f\left( y\right)\Phi\left( y\right)
\right)\right.  \notag \\
&-\left. f\left( x\right) \left( \sum_{y\sim x} a_{\lbrace x,y\rbrace} \left[
\Phi\left( x\right)-\Phi\left(y \right) \right] \right) f\left(
x\right)\Phi\left( x\right) \right]  \notag \\
&=\sum_{x\in V}f\left( x\right)\Phi\left( x\right) \sum_{y\sim x} a_{\lbrace
x,y\rbrace} \Phi\left( y\right)\left[ f\left( x\right)-f\left(y \right) %
\right]  \notag \\
&=\sum_{x\in V}\Phi^{2}\left( x\right)f\left( x\right) \dfrac{1}{%
\Phi^2\left( x\right)} \sum_{y\sim x} a_{\lbrace x,y\rbrace}\Phi\left(
x\right)\Phi\left( y\right) \left[ f\left( x\right)-f\left(y \right) \right]
\notag
\end{align}

En posant $$c_{\left\lbrace x,y\right\rbrace}= a_{\lbrace
x,y\rbrace}\Phi\left( x\right)\Phi\left( y\right)~,$$
 il r\'esulte que:
\begin{equation*}
\left\langle Pg,g\right\rangle_{l^2}= \left\langle
\Delta_{\omega,c}f,f\right\rangle_{l^2_{\omega}}
\end{equation*}
D'o\`u
\begin{equation*}
P=U^{-1}\Delta_{\omega,c}U~,
\end{equation*}
avec $U:l^2 \longrightarrow l^2_{\omega}$ d\'efinie par $U\left( g\right)=%
\dfrac{g}{\Phi}$~.\newline
Ainsi $P$ est unitairement \'equivalent au Laplacien $\Delta_{\omega,c}$,
avec $\omega \equiv \Phi$ et $c$ donn\'e par $c_{\left\lbrace x,y
\right\rbrace}=a_{\left\lbrace x,y \right\rbrace}\Phi(x)\Phi(y)$~.
\end{demorefunit}

\section{Cas des graphes m\'etriquement complets}\label{complesa}
Nous adaptons au cas des graphes la m\'ethode de G. et I. Nenciu \cite{Nen}
utilisant une estimation d'Agmon donn\'ee dans le lemme technique suivant
utile pour la preuve du th\'eor\`eme \ref{compl}~.

\begin{lemm}
\label{estim} Soient $H =\Delta_{1,a}+W$ un op\'erateur de Schr\"odinger sur
$G$, $\lambda$ un nombre r\'eel et $v\in C\left(V \right)$~.\newline
On suppose que $v$ est une solution de l'\'equation: $\left(
H-\lambda\right)\left( v\right) =0 $~. Alors pour tout $f\in C_{0}\left(
V\right)$~, on a:
\begin{align}
\left\langle fv,\left( H-\lambda\right) \left( fv\right) \right\rangle_{l^{2}}
&=\sum_{\{x,y \} \in E}a_{\lbrace x,y\rbrace}v\left(x
\right)v\left( y\right) \left[ f\left(x \right)-f\left( y\right) \right]^{2}
\notag \\
&=\dfrac{1}{2}\sum_{x\in V}v\left( x\right) \sum _{y\sim x} a_{\lbrace
x,y\rbrace}v\left( y\right) \left[ f\left(x \right)-f\left( y\right) \right]%
^{2}~.  \notag
\end{align}
\end{lemm}

\begin{demo}
On suppose que: $\left( H-\lambda\right)\left( v\right)=0$~, c'est \`a dire
qu'on a pour tout $x\in V$~,
\begin{equation*}
\sum _{y\sim x} a_{\lbrace x,y\rbrace}\left( v\left( x\right) -v\left(
y\right)\right) + W\left( x\right)v\left( x\right)=\lambda v\left( x\right)
\end{equation*}
Calculons $S=\left\langle fv,\left( H-\lambda\right) \left( fv\right)
\right\rangle_{l^{2}}$
\begin{align}
S &= \sum_{x\in V}f\left( x\right)v\left( x\right) \left[\left( H-\lambda
\right)\left( fv\right)\right]\left(x \right)  \notag \\
&=\sum_{x\in V}f\left( x\right)v\left( x\right) W\left(x\right)f\left( x\right)v\left( x\right) -
\lambda f\left( x\right)v\left(x\right) \notag \\
&+ \sum_{x\in V} \sum _{y\sim x} a_{\lbrace x,y\rbrace} \left[f\left(
x\right)v\left( x\right)-f\left( y\right)v\left( y\right) \right]~.
\notag
\end{align}
Or par l'hypoth\`ese faite sur $v$~, on a:
\begin{equation*}
\lambda f\left( x\right)v\left( x\right)-W\left( x\right)f\left(
x\right)v\left( x\right) = \sum _{y\sim x} a_{\lbrace x,y\rbrace} f\left(
x\right)\left[ v\left( x\right)-v\left( y\right)\right]~.
\end{equation*}
On obtient alors, en rempla\c{c}ant dans l'expression pr\'ec\'edente:
\begin{align}
S &=\sum _{x\in V}f\left( x\right)v\left( x\right)\sum _{y\sim x} a_{\lbrace
x,y\rbrace}v\left( y\right) \left[ f\left(x \right)-f\left( y\right) \right]
\notag \\
&=\sum _{x\in V}\sum _{y\sim x} a_{\lbrace x,y\rbrace}v\left(x \right)
v\left( y\right) \left[ f^{2}\left(x \right)-f\left( x\right) f\left(
y\right) \right] \notag
\end{align}
Comme $a_{\left\lbrace x,y\right\rbrace }=a_{\left\lbrace y,x\right\rbrace }~,$
 l'expression devient:
\begin{equation*}
S =\sum_{\{x,y \} \in E} a_{\lbrace x,y\rbrace}v\left(x
\right)v\left( y\right) \left[ f^{2}\left( x\right) -f\left( x\right)
f\left( y\right) + f^{2}\left( y\right) -f\left( x\right) f\left( y\right)%
\right]
\end{equation*}
D'o\`u finalement:
\begin{align}
\left\langle fv,\left( H-\lambda\right) \left( fv\right) \right\rangle_{l^{2}}
&=\sum_{\{x,y \} \in E}a_{\lbrace x,y\rbrace}v\left(x
\right)v\left( y\right) \left[ f\left(x \right)-f\left( y\right) \right]^{2}
\notag \\
&=\dfrac{1}{2}\sum_{x\in V}v\left( x\right) \sum _{y\sim x} a_{\lbrace
x,y\rbrace}v\left( y\right) \left[ f\left(x \right)-f\left( y\right) \right]%
^{2}~.  \notag
\end{align}
\end{demo}\newline

\begin{defi}
Un graphe $G$ est dit \`a valence born\'ee, s'il existe un entier $N$ tel que pour
tout $x\in V$ on ait: $\sharp \left\lbrace y\in V;y\sim x\right\rbrace\leq N~.$
\end{defi}

\begin{defi}
Soit $a$ une fonction strictement positive sur les ar\^etes de $G$~. On
d\'efinit la distance pond\'er\'ee par $a$ sur $G$~, qu'on note $\delta_{a}$~,
par:
\begin{equation*}
\delta_{a} \left( x,y\right)=\min _{\gamma\in \Gamma_{x,y}}L\left( \gamma
\right)
\end{equation*}
o\`u $\Gamma_{x,y} $ est l'ensemble de tous les chemins d'ar\^etes $%
\gamma:x_{1}=x,x_2,$...$x_{n}=y$~, reliant $x$ \`a $y$~; et $L\left( \gamma
\right)=\underset{1\leq i\leq n}{\sum} \dfrac{1}{\sqrt{a_{x_{i}x_{i+1}}}}$
la longueur du chemin d'ar\^etes $\gamma$~.
\end{defi}

\begin{theo}
\label{compl}  Soit $H=\Delta_{1,a}+W$ un op\'erateur de Schr\"odinger sur
un graphe infini $G$ \`a valence born\'ee et telle que sa
m\'etrique d\'efinie par la distance $\delta_{a}$ est compl\`ete. On suppose
qu'il existe un r\'eel $k$ telle que $$\left\langle Hg,g
\right\rangle_{l^{2}} \geq k\Vert g\Vert^{2}_{l^{2}}$$~
pour tout $g\in C_{0 }\left( V\right)$.
Alors l'op\'erateur $H$~, avec comme domaine $C_{0 }\left( V\right)
$, est essentiellement auto-adjoint.
\end{theo}

\begin{demo}
Soit $\lambda \ < k-1$~, on va montrer que si $v\in l^{2}\left( V\right)$ et
v\'erifie l'\'equation $Hv=\lambda v$~, alors $v$ est identiquement nulle.%
\newline
On fixe $R\ >0$ et un sommet $x_{0}$ qu'on prend comme origine. On note:
$$B_{R}=\left\lbrace x\in V;\delta_{a} \left( x_{0},x\right )\leq
R\right\rbrace $$ la boule de centre $x_{0}$ et de rayon $R$ pour la distance
$\delta_{a}$.\newline
On consid\`ere la fonction $f$ d\'efinie sur $V$ par
\begin{equation*}
f\left( x\right) =\min\left( 1,\delta_{a}\left(x, V\setminus
B_{R+1}\right)\right)~.
\end{equation*}
On a ainsi:
\begin{equation*}
f\upharpoonright B_{R}\equiv 1,\quad f\upharpoonright V\setminus
B_{R+1}\equiv0~,  \quad f\left( B_{R+1}\setminus B_{R} \right)\subseteq \left[
0,1\right]
\end{equation*}
Le support de $f$ est inclus dans $B_{R+1}$ qui est fini du fait que la
m\'etrique associ\'ee \`a $\delta_{a}$ est compl\`ete.\newline
Par l'hypoth\`ese faite sur $H$ et vu que $fv$ est \`a support fini dans $%
B_{R+1}$~, on obtient les minorations suivantes:
\begin{equation*}
\left\langle fv,\left( H-\lambda\right) \left( fv\right) \right\rangle_{l^{2}} \geq
\left( k-\lambda \right) \sum_{x\in B_{R+1}}\left( fv\right) ^2\left( x\right)
\geq\sum _{x\in B_{R}} v^2\left( x\right)
\end{equation*}
D'autre part, en utilisant le lemme \ref{estim}~, on obtient:
\begin{align}
\left\langle fv,\left( H-\lambda \right) \left( fv\right) \right\rangle_{l^{2}}&=
\dfrac{1}{2}\sum _{x\in V}\sum_{y\sim x}a_{\left\lbrace x,y\right\rbrace}
v\left( x\right)v\left(y \right) \left[f\left( x\right)-f\left( y\right) %
\right]^2  \notag \\
&\leq \dfrac{1}{2}\sum _{x\in V}\sum_{y\sim x}a_{\left\lbrace
x,y\right\rbrace} v^{2}\left( x\right) \left[f\left( x\right)-f\left(
y\right) \right]^2  \notag
\end{align}
en remarquant que $a_{\left\lbrace x,y\right\rbrace}=a_{\left\lbrace
y,x\right\rbrace}$ et que:
\begin{equation*}
v\left( x\right)v\left(y \right) \leq\dfrac{1}{2}\left(v^{2}\left(
x\right)+v^{2}\left( y\right) \right)~.
\end{equation*}
Comme chacune des restrictions de $f$ \`a $B_{R}$ et \`a $V\setminus B_{R+1}$
sont des fonctions constantes, alors l'in\'egalit\'e pr\'ec\'edente devient:
\begin{align}
\left\langle fv,\left( H-\lambda\right) \left( fv\right) \right\rangle_{l^{2}}&\leq
\dfrac{1}{2}\sum _{x\in B_{R+1}\setminus B_{R}}\sum_{y\sim x}a_{\left\lbrace
x,y\right\rbrace} v^2\left( x\right) \left[f\left( x\right)-f\left( y\right) %
\right]^2  \notag \\
&\leq\dfrac{1}{2}\sum _{x\in B_{R+1}\setminus B_{R}}\sum_{y\sim
x}a_{\left\lbrace x,y\right\rbrace} v^{2}\left( x\right)\left(
\delta_{a}\left( x,y\right) \right)^{2}  \notag
\end{align}
La derni\`ere in\'egalit\'e est obtenue du fait que $f$ est $1-$%
Lipschitzienne puisque qu'elle est minimum de deux fonctions $1-$%
Lipschitziennes.\newline
Etant donn\'e que la distance $\delta_{a}$ satisfait l'in\'egalit\'e:
\begin{equation*}
\delta_{a}\left( x,y\right)\leq\frac{1}{\sqrt{a_{\left\lbrace
x,y\right\rbrace}}}
\end{equation*}
si $\left\lbrace x,y\right\rbrace$ est une ar\^ete, et que la valence de $G$
est uniform\'ement born\'ee par $N$~, on trouve:
\begin{equation*}
\left\langle fv,\left( H-\lambda\right) \left( fv\right) \right\rangle_{l^{2}} \leq%
\dfrac{1}{2}N\sum_{x\in B_{R+1}\setminus B_{R}}v^2\left( x\right)
\end{equation*}
Ainsi, pour tout $R\ >0$~, on obtient:
\begin{equation*}
\sum _{x\in B_{R}}v^2\left( x\right)\leq\left\langle fv,\left(
H-\lambda\right) \left( fv\right) \right\rangle_{l^{2}} \leq\dfrac{1}{2}N\sum_{x\in
B_{R+1}\setminus B_{R}}v^2\left( x\right)
\end{equation*}
Et puisque $v\in l^2\left( V\right)$~, en faisant tendre $R\rightarrow\infty$~,
 on obtient:
\begin{equation*}
\lim_{R\rightarrow\infty}\sum _{x\in B_{R+1}\setminus B_{R}}v^2\left(
x\right) =0
\end{equation*}
et par suite $\Vert v\Vert^2_{l^2} =0$~.
\end{demo}

\begin{theo}\label{lapcompl}
Soit $G$ un graphe infini \`a valence uniform\'ement born\'ee et $%
\Delta_{\omega,c}$ un Laplacien sur $G$. On suppose que la m\'etrique
associ\'ee \`a la distance $\delta_{a}$ est compl\`ete, o\`u $a$ est la
fonction donn\'ee par $$a_{{\left\lbrace x,y\right\rbrace }}= \dfrac{%
c_{\left\lbrace x,y\right\rbrace }}{\omega_{x}\omega_{y}}~.$$\newline
Alors l'op\'erateur $\Delta _{\omega,c~,}$ avec domaine $C_0(V)$, est
essentiellement auto-adjoint.
\end{theo}

\begin{demo}
Par le lemme \ref{lap-unit} l'op\'erateur $\Delta _{\omega,c}$ est
unitairement \'equivalent \`a l'op\'erateur de Schr\"odinger $%
H=\Delta_{1,a}+W$~, o\`u
\begin{equation*}
a_{\left\lbrace x,y\right\rbrace }= \dfrac{c_{\left\lbrace x,y\right\rbrace }%
}{\omega_{x}\omega_{y}}
\end{equation*}
et
\begin{equation*}
W\left( x\right)= \dfrac{1}{\omega_{x}^{2}}~\sum_{y\sim x}
c_{\lbrace x,y\rbrace} \left( 1-\dfrac{\omega_{x}}{\omega_{y}}\right)
\end{equation*}
Et on applique le th\'eor\`eme \ref{compl} \`a l'op\'erateur $H$ qui
v\'erifie les hypoth\`eses, puisque $\Delta_{\omega,c}$ est positif et qu'on
a $$\left\langle Hg,g \right\rangle_{l^2}= \left\langle \Delta_{\omega,c}%
\dfrac{g}{\omega} ,\dfrac{g}{\omega} \right\rangle_{l_{\omega}^2}$$ dans
la preuve du th\'eor\`eme \ref{unit} .
\end{demo}

\begin{rem}\label{non-min}
Le th\'eor\`eme \ref{compl} n'est pas un cas particulier du th\'eor\`eme
\ref{op.eaa}. En effet dans le th\'eor\`eme \ref{compl} le potentiel $W$ n'est
pas n\'ecessairement minor\'e. En effet, prenons par exemple
$G$ tel que $V= {\mathbb{N}} \setminus\ \left\lbrace
0,1\right\rbrace $ et $n\sim n+1$ pour tout $n$~. Supposons que $G$ est
pond\'er\'e par \; $\omega_n=\dfrac{1}{n\ \log n}$ sur $V$ et connect\'e par\;  $%
c_n=1$ sur $E$~. La distance $\delta_{a}$ est donn\'ee par:
\begin{equation*}
\delta_{a}( n_0,n )=\sum_{n_0\leq k\leq n} \dfrac{1}{\sqrt{a_{k,k+1}}}
\end{equation*}
or
\begin{equation*}
\dfrac{1}{\sqrt{a_{k,k+1}}} = \dfrac{1}{\sqrt{k(k+1)\log k\log(k+1)}}%
\underset{\infty}{\sim}\dfrac{1}{k\log k}~,
\end{equation*}
donc   $\delta_{a}\ ( n_0,n\ )\underset{n\rightarrow\infty}{%
\longrightarrow}\infty$~, et la m\'etrique associ\'ee est compl\`ete.
\newline
De plus, en posant $H=\Delta_{1,a}+W$~, nous avons: $\left\langle Hg,g
\right\rangle_{l^2}= \left\langle \Delta_{\omega,c}\dfrac{g}{\omega},%
\dfrac{g}{\omega} \right\rangle_{l_{\omega}^2}\geq0$, pour $g\in C_0(V)~.$
Alors que le potentiel $W$ n'est pas minor\'e, puisqu'apr\`es calculs, nous
obtenons:
\begin{equation*}
W\left( n\right)=2n^2\log^2n-n\log n\left[(n+1)\log(n+1)+ (n-1)\log(n-1)%
\right]\underset{\infty}{\sim}-\log n
\end{equation*}
qui tend vers $-\infty$~.
\end{rem}

\begin{rem}
Dans l'exemple de la remarque \ref{non-min}~, le choix du poids en fonction de $log$ est
d\'eterminant. En effet, en prenant des fonctions puissances, on ne peut
pas avoir \`a la fois la m\'etrique $\delta_{a}\ $ compl\`ete et
le potentiel $W$ non minor\'e. Prenons par exemple $G$ tel que $V= {\mathbb{N}} \setminus\ \left\lbrace
0,1\right\rbrace $ et $n\sim n+1$ pour tout $n$~. Supposons que $G$ soit
pond\'er\'e par le poids\; $\omega_n=\dfrac{1}{n^{\alpha}}$ sur $V$ et la conductance
$c_n=\dfrac{1}{n^{\beta}}$ sur $E$~. La distance $\delta_{a}$ est donn\'ee par:

\begin{equation*}
\delta_{a}( n_0,n )=\sum_{n_0\leq k\leq n}\dfrac{k^{\beta}}{(k^{\alpha}(k+1)^{\alpha})^{\frac{1}{2}}}~.
\end{equation*}
Et un simple calcul montre que
$$\delta_{a}~ \textnormal{est compl\`ete \ssi}\ \alpha-\dfrac{1}{2}~\beta \leq 1~.$$
Pour ce qui est du potentiel, le calcul donne:
$$W_{n}\sim -\alpha(\alpha -\beta-1)n^{2\alpha-\beta-2}~.$$
Ainsi pour $\alpha-\dfrac{1}{2}~\beta \leq 1~,$ $W_{n}$
est constant ou bien tend vers $0~,$ donc il est minor\'e.
\end{rem}

\begin{rem}
La condition que la m\'etrique $\delta_{a}$ soit compl\`ete n'est pas
une condition n\'ecessaire pour que le Laplacien $\Delta_{\omega,c}$ ou
que l'op\'erateur de Schr\"{o}dinger $\Delta_{1,a}$ soit
\ESA. En effet,
soit $G$ tel que $V= {\mathbb{N}} \setminus\ \left\lbrace
0,1\right\rbrace $ et $n\sim n+1$ pour tout $n$~.  Pour tout poids
$a$ sur $V$ le laplacien $\Delta_{1,a}$
est \ESA par le th\'eor\`eme \ref{lapeaa}~. Alors que la m\'etrique $\delta_{a}$
n'est pas n\'ecessairement compl\`ete: par exemple lorsque
$a_{n}=\dfrac{1}{n^{2+\varepsilon}}$~, pour $\varepsilon>0~.$
\newline
Dans l'article [C-To-Tr-1]
on donnera des conditions de croissance du potentiel
assurant qu'un op\'erateur de Schr\"{o}dinger
est essentiellement auto-adjoint dans le cas
des graphes m\'etriquement non complets.
\end{rem}
\begin{remer}
Ce travail a \'et\'e r\'ealis\'e gr\^ace aux financements par l'unit\'e de recherches (05/UR/15-02) 
"Math\'ematiques et Applications" et par le Département de Math\'ematiques de la Facult\'e des
Sciences de Bizerte (Tunisie) pour mes multiples courts s\'ejours de recherches  au Laboratoire de l'Institut
Fourier (Grenoble-France).
\newline
Je tiens \`a exprimer ma gratitude au Professeur Yves Colin de Verdi\`ere pour sa
collaboration et pour avoir formul\'e l'id\'ee de ce travail. Je remercie le
Professeur Fran\c{c}oise Truc pour l'int\'er\^eret qu'elle a port\'e \`a ce travail et pour les
 multiples \'echanges fructueux.
\end{remer}

\end{document}